\documentclass[reqno]{amsart}
\usepackage{amsmath,amssymb,amsfonts,latexsym,txfonts,pxfonts,wasysym}
\usepackage{amssymb,longtable}
\usepackage{hyperref}
\usepackage{lscape,graphicx}
\begin{document}
\title[ Two Results due to Ramanujan ]
{Remark On Two Results due to Ramanujan on Hypergeometric Series}

\author[M. A. Rakha, A. K. Ibrahim, A. K. Rathie]
{Medhat A. Rakha, Adel K. Ibrahim, Arjun K. Rathie}  

\address{Medhat A. Rakha \newline
 Department of Mathematics and Statistics,
 College of Science,
 Sultan Qaboos University,
 P.O.Box 36 - Al-Khoud 123, 
 Muscat - Sultanate of Oman}
\email{medhat@squ.edu.om}

\address{Adel K. Ibrahim \newline
 Mathematics Department, 
 College of Science,
 Jazan University University, Jazan, Saudi Arabia}  
 \email{dradlkhalil@yahoo.com}

\address{Arjun K. Rathie \newline
 Department of Mathematics, 
 School of Mathematical and Physical Sciences,
 Central University of Kerala, Riverside Transit Campus,
 Padennakkad P.O. Nileshwar, Kasaragod - 671 328, Kerala - INDIA}
 \email{akrathie@gmail.com}


\begin{abstract}
During the course of verifying the results of Ramanujan on hypergeometric
series, Berndt in his notebooks, Part II mentioned corrected forms of two of
the Ramanujan's results.

The aim of this short research note is to point out that one of the results
obtained by Ramanujan is correct (and not of Berndt's result) and the second
result (which is slightly differ from Ramanujan's result and Berndt's result)
is given here in corrected form.

\end{abstract}

\maketitle
\numberwithin{equation}{section}
\newtheorem{theorem}{Theorem}[section]
\newtheorem{lemma}[theorem]{Lemma}
\newtheorem{proposition}[theorem]{Proposition}
\newtheorem{corollary}[theorem]{Corollary}
\newtheorem*{remark}{Remark}

\section{Introduction}
By defining
\[
\mu=\frac{\Gamma\left(\frac{1}{2}\right)}{\Gamma^{2}\left(\frac{3} {4}\right)}~
\text{~~and~~}
\eta=\frac{\Gamma^{2}\left(  \frac{3}{4}\right)  }{\Gamma^{3}\left(  \frac
{1}{2}\right)  }%
\]
we start with the following two results given by Ramanujan \cite[p. 141 and p.
142]{4}%
\begin{align}
&  _{2}F_{1}\left[
\begin{array}
[c]{ccc}%
\frac{1}{2}, & \frac{1}{2}; & \\
&  & \frac{1}{2}+\frac{x}{1+x^{2}}\\
1; &  &
\end{array}
\right]  \nonumber\\
&  =\mu\sqrt{1+x^{2}}\,_{2}F_{1}\left[
\begin{array}
[c]{ccc}%
\frac{1}{4}, & \frac{1}{2}; & \\
&  & x^{4}\\
\frac{3}{4}; &  &
\end{array}
\right]  +\eta x\sqrt{1+x^{2}}\,_{2}F_{1}\left[
\begin{array}
[c]{ccc}%
\frac{3}{4}, & \frac{1}{2}; & \\
&  & x^{4}\\
\frac{5}{4}; &  &
\end{array}
\right]  \label{1}%
\end{align}
and%
\begin{align}
&  \sqrt{1-x^{2}}\,_{2}F_{1}\left[
\begin{array}
[c]{ccc}%
\frac{1}{2}, & \frac{1}{2}; & \\
&  & \frac{1}{2}+\frac{x}{1+x^{2}}\\
1; &  &
\end{array}
\right]  \nonumber\\
&  =\mu\,_{2}F_{1}\left[
\begin{array}
[c]{ccc}%
\frac{1}{2}, & \frac{1}{2}; & \\
&  & \frac{x^{4}}{x^{4}-1}\\
\frac{3}{4}; &  &
\end{array}
\right]  +2\eta x\,_{2}F_{1}\left[
\begin{array}
[c]{ccc}%
\frac{1}{2}, & \frac{1}{2}; & \\
&  & \frac{x^{4}}{x^{4}-1}\\
\frac{5}{4}; &  &
\end{array}
\right]  .\label{2}%
\end{align}

Berndt \cite{1}, in his Ramanujan's notebooks, Part II pointed out that these
results contain some errors and should respectively be read as

\begin{itemize}
\item Entry 34 (ii) \cite[eq. (34.3). p.97]{1}%
\begin{align}
&  _{2}F_{1}\left[
\begin{array}
[c]{ccc}%
\frac{1}{2}, & \frac{1}{2}; & \\
&  & \frac{1}{2}+\frac{x}{1+x^{2}}\\
1; &  &
\end{array}
\right] \nonumber\\
&  =\mu\sqrt{1+x^{2}}\,_{2}F_{1}\left[
\begin{array}
[c]{ccc}%
\frac{1}{4}, & \frac{1}{2}; & \\
&  & x^{4}\\
\frac{3}{4}; &  &
\end{array}
\right]  +\eta x\left(  1+x^{2}\right)  ^{\frac{3}{2}}\,_{2}F_{1}\left[
\begin{array}
[c]{ccc}%
\frac{3}{4}, & \frac{1}{2}; & \\
&  & x^{4}\\
\frac{5}{4}; &  &
\end{array}
\right]  \label{3}%
\end{align}
and

\item Example (ii) \cite[p. 99]{1}%
\begin{align}
&  \sqrt{1-x^{2}}\,_{2}F_{1}\left[
\begin{array}
[c]{ccc}%
\frac{1}{2}, & \frac{1}{2}; & \\
&  & \frac{1}{2}+\frac{x}{1+x^{2}}\\
1; &  &
\end{array}
\right]  \nonumber\\
&  =\mu\,_{2}F_{1}\left[
\begin{array}
[c]{ccc}%
\frac{1}{2}, & \frac{1}{2}; & \\
&  & \frac{x^{4}}{x^{4}-1}\\
\frac{3}{4}; &  &
\end{array}
\right]  +\eta x\,(1+x^{2})_{2}F_{1}\left[
\begin{array}
[c]{ccc}%
\frac{1}{2}, & \frac{1}{2}; & \\
&  & \frac{x^{4}}{x^{4}-1}\\
\frac{5}{4}; &  &
\end{array}
\right]  .\label{4}%
\end{align}

\end{itemize}

In this short research note, we shall show that Ramanujan's result (\ref{2})
is \textbf{correct  } (and not (\ref{4}) obtained by Brendt) and the corrected
form of the Ramanujan's result (\ref{1}) (and of Brendt's result (\ref{3}))
should be read as%
\begin{align}
&  _{2}F_{1}\left[
\begin{array}
[c]{ccc}%
\frac{1}{2}, & \frac{1}{2}; & \\
&  & \frac{1}{2}+\frac{x}{1+x^{2}}\\
1; &  &
\end{array}
\right]  \nonumber\\
&  =\mu\sqrt{1+x^{2}}\,_{2}F_{1}\left[
\begin{array}
[c]{ccc}%
\frac{1}{2}, & \frac{1}{4}; & \\
&  & x^{4}\\
\frac{3}{4}; &  &
\end{array}
\right]  +2\eta x\sqrt{1+x^{2}}\,_{2}F_{1}\left[
\begin{array}
[c]{ccc}%
\frac{3}{4}, & \frac{1}{2}; & \\
&  & x^{4}\\
\frac{5}{4}; &  &
\end{array}
\right]  .\label{5}%
\end{align}

In order to derive these results, we start with the following known Kummer's
Formula \cite[p. 64]{1}%
\begin{align}
&  _{2}F_{1}\left[
\begin{array}
[c]{ccc}%
a, & b; & \\
&  & \frac{1+x}{2}\\
\frac{1}{2}(a+b+1); &  &
\end{array}
\right]  \nonumber\\
&  =\frac{\Gamma\left(  \frac{1}{2}\right)  \Gamma\left(  \frac{1}{2}%
a+\frac{1}{2}b+\frac{1}{2}\right)  }{\Gamma\left(  \frac{1}{2}a+\frac{1}%
{2}\right)  \Gamma\left(  \frac{1}{2}b+\frac{1}{2}\right)  }\,_{2}F_{1}\left[
\begin{array}
[c]{ccc}%
\frac{1}{2}a, & \frac{1}{2}b; & \\
&  & x^{2}\\
\frac{1}{2}; &  &
\end{array}
\right]  \nonumber\\
&  +\frac{2x\Gamma\left(  \frac{1}{2}\right)  \Gamma\left(  \frac{1}{2}%
a+\frac{1}{2}b+\frac{1}{2}\right)  }{\Gamma\left(  \frac{1}{2}a\right)
\Gamma\left(  \frac{1}{2}b\right)  }\,_{2}F_{1}\left[
\begin{array}
[c]{ccc}%
\frac{1}{2}a+\frac{1}{2}, & \frac{1}{2}b+\frac{1}{2}; & \\
&  & x^{2}\\
\frac{3}{2}; &  &
\end{array}
\right]  .\label{6}%
\end{align}

It is not out of place to mention here that result (\ref{6}) was also
independently rediscovered by Ramanujan \cite[Entry 21, p. 64]{1}.

Further, in (\ref{6}), if we take $a=b=\frac{1}{2}$, we get the following
result due to Ramanujan \cite[Entry 34 (i), p. 96]{1}%
\begin{align}
& _{2}F_{1}\left[
\begin{array}
[c]{ccc}%
\frac{1}{2}, & \frac{1}{2}; & \\
&  & \frac{1+x}{2}\\
1; &  &
\end{array}
\right]  \nonumber\\
& =\mu\,_{2}F_{1}\left[
\begin{array}
[c]{ccc}%
\frac{1}{4}, & \frac{1}{4}; & \\
&  & x^{2}\\
\frac{1}{2}; &  &
\end{array}
\right]  +\eta x\,_{2}F_{1}\left[
\begin{array}
[c]{ccc}%
\frac{3}{4}, & \frac{3}{4}; & \\
&  & x^{2}\\
\frac{3}{2}; &  &
\end{array}
\right]  .\label{7}%
\end{align}

Also, we mention here another result due to Kummer, which was also
rediscovered by Ramanujan \cite[Entry 3, p. 50]{1} viz.%

\begin{equation}
_{2}F_{1}\left[
\begin{array}
[c]{ccc}%
r, & m; & \\
&  & \frac{4x}{(1+x)^{2}}\\
2m; &  &
\end{array}
\right]  =(1+x)^{2r}\,_{2}F_{1}\left[
\begin{array}
[c]{ccc}%
r, & r-m+\frac{1}{2}; & \\
&  & x^{2}\\
m+\frac{1}{2}; &  &
\end{array}
\right]  .\label{8}%
\end{equation}

In (\ref{8}), if we replace $x$ by $y^{2}$ and take (i) $\ r=m=\frac{1}{4}$
and (ii) $r=m=\frac{3}{4}$, we respectively get%

\begin{equation}
_{2}F_{1}\left[
\begin{array}
[c]{ccc}%
\frac{1}{4}, & \frac{1}{4}; & \\
&  & \frac{4y^{2}}{(1+y^{2})^{2}}\\
\frac{1}{2}; &  &
\end{array}
\right]  =\sqrt{1+y^{2}}\,_{2}F_{1}\left[
\begin{array}
[c]{ccc}%
\frac{1}{4}, & \frac{1}{2}; & \\
&  & y^{4}\\
\frac{3}{4}; &  &
\end{array}
\right]  ,\label{9}%
\end{equation}
and%

\begin{equation}
_{2}F_{1}\left[
\begin{array}
[c]{ccc}%
\frac{3}{4}, & \frac{3}{4}; & \\
&  & \frac{4y^{2}}{(1+y^{2})^{2}}\\
\frac{3}{2}; &  &
\end{array}
\right]  =\left(  1+y^{2}\right)  ^{\frac{3}{2}}\,_{2}F_{1}\left[
\begin{array}
[c]{ccc}%
\frac{3}{4}, & \frac{1}{2}; & \\
&  & y^{4}\\
\frac{5}{4}; &  &
\end{array}
\right]  .\label{10}%
\end{equation}

Also, the well known Euler's first transformation \cite[Theorem 20, p. 60]{3}%

\begin{equation}
_{2}F_{1}\left[
\begin{array}
[c]{ccc}%
a, & b; & \\
&  & x\\
c; &  &
\end{array}
\right]  =\left(  1-x\right)  ^{-a}\,_{2}F_{1}\left[
\begin{array}
[c]{ccc}%
a, & c-b; & \\
&  & -\frac{x}{1-x}\\
c; &  &
\end{array}
\right]  .\label{11}%
\end{equation}

Now with these results, we are ready to establish Ramanujan's results:

\begin{proof}
In (\ref{7}), replacing $x$ by $\frac{2y}{1+y^{2}}$, we have
\begin{align}
& _{2}F_{1}\left[
\begin{array}
[c]{ccc}%
\frac{1}{2}, & \frac{1}{2}; & \\
&  & \frac{1}{2}+\frac{y}{1+y^{2}}\\
1; &  &
\end{array}
\right]  \nonumber\\
& =\mu\,_{2}F_{1}\left[
\begin{array}
[c]{ccc}%
\frac{1}{4}, & \frac{1}{4}; & \\
&  & \frac{4y^{2}}{\left(  1+y^{2}\right)  ^{2}}\\
\frac{1}{2}; &  &
\end{array}
\right]  +\eta\frac{2y}{1+y^{2}}\,_{2}F_{1}\left[
\begin{array}
[c]{ccc}%
\frac{3}{4}, & \frac{3}{4}; & \\
&  & \frac{4y^{2}}{\left(  1+y^{2}\right)  ^{2}}\\
\frac{3}{2}; &  &
\end{array}
\right]  .\label{12}%
\end{align}

Upon using the results given by (\ref{9}) and (\ref{10}) on the right-hand
side of (\ref{12}), we get the corrected form of the result (\ref{1}) due to
Ramanujan and (\ref{3}) due Brendt as
\begin{align}
& _{2}F_{1}\left[
\begin{array}
[c]{ccc}%
\frac{1}{2}, & \frac{1}{2}; & \\
&  & \frac{1}{2}+\frac{y}{1+y^{2}}\\
1; &  &
\end{array}
\right]  \nonumber\\
& =\mu\sqrt{1+y^{2}}\,_{2}F_{1}\left[
\begin{array}
[c]{ccc}%
\frac{1}{4}, & \frac{1}{2}; & \\
&  & y^{4}\\
\frac{3}{4}; &  &
\end{array}
\right]  +2\eta\,y\sqrt{1+y^{2}}\,_{2}F_{1}\left[
\begin{array}
[c]{ccc}%
\frac{3}{4}, & \frac{1}{2}; & \\
&  & y^{4}\\
\frac{5}{4}; &  &
\end{array}
\right]  .\label{13}%
\end{align}

Further, upon applying Euler's first transformation formula (\ref{11}) in the
right-hand side of (\ref{13}), we at once get (\ref{2}), which is the
Ramanujan result.

We conclude the note by remarking that the results (\ref{2}) and (\ref{13})
have been checked and verified numerically through \textit{Mathematica}.
\end{proof}

\subsection*{Acknowledgments}
The work of this research was  supported by the research grant (05/4/33) funded by Jazan University - Jazan, Saudi Arabia.

\end{document}